\documentclass[12pt]{article}
\usepackage{amssymb,amsthm,amsmath}

\theoremstyle{plain}
\newtheorem{Theorem}{Theorem}[section]
\newtheorem{Lemma}[Theorem]{Lemma}

\newtheorem{Proposition}[Theorem]{Proposition}
\newtheorem{Corollary}[Theorem]{Corollary}
\newtheorem{Conjecture}[Theorem]{Conjecture}
\newtheorem{Question}[Theorem]{Question}

\newtheorem*{Theorem*}{Theorem}
\newtheorem*{Lemma*}{Lemma}

\theoremstyle{definition}
\newtheorem*{Definition*}{Definition}
\newtheorem*{Remark*}{Remark}
\newtheorem*{Example*}{Example}
\newtheorem{Definition}[Theorem]{Definition}
\newtheorem{Remark}[Theorem]{Remark}

\newcommand{\isom}{\cong}

\newcommand{\gen}[1]{\mathopen{<}#1\mathclose{>}}
\newcommand{\abs}[1]{\lvert#1\rvert}

\DeclareMathOperator{\Gal}{Gal}

\DeclareMathOperator{\ind}{ind}

\newcommand{\CCC}{\mathbb C}
\newcommand{\NNN}{\mathbb N}
\newcommand{\QQQ}{\mathbb Q}

\newcommand{\ZZZ}{\mathbb Z}
\newcommand{\PP}{\mathbb P^1}           

\newcommand{\hk}{\hat k}
\newcommand{\bQ}{\overline\QQQ}

\newcommand{\ok}{{\mathcal O}_k}        

\newcommand{\frakp}{\mathfrak p}
\newcommand{\frakP}{\mathfrak P}

\DeclareMathOperator{\GL}{GL}

\DeclareMathOperator{\PSL}{PSL}

\DeclareMathOperator{\AGL}{AGL}

\newcommand{\M}[2]{\mathrm{M}_{#1#2}}             

\DeclareMathOperator{\alt}{\mathcal{A}}
\DeclareMathOperator{\sym}{\mathcal{S}}
\DeclareMathOperator{\Red}{Red}

\newcommand{\bk}{\bar k}

\newcommand{\Si}{S}
\newcommand{\alphlab}{
  \renewcommand{\theenumi}{(\alph{enumi})}
  \renewcommand{\labelenumi}{\theenumi}
  \renewcommand{\theenumii}{(\roman{enumii})}
  \renewcommand{\labelenumii}{\theenumii}
}
\allowdisplaybreaks[1]

\begin{document}
\title{Finiteness Results for Hilbert's Irreducibility Theorem}
\author{Peter M\"uller}
\maketitle

\alphlab

\begin{abstract}
  Let $k$ be a number field, $\ok$ its ring of integers, and
  $f(t,X)\in k(t)[X]$ be an irreducible polynomial. Hilbert's
  irreducibility theorem gives infinitely many integral
  specializations $t\mapsto\bar t\in\ok$ such that $f(\bar t,X)$ is
  still irreducible. In this paper we study the set $\Red_f(\ok)$ of
  those $\bar t\in\ok$ with $f(\bar t,X)$ reducible. We show that
  $\Red_f(\ok)$ is a finite set under rather weak assumptions. In
  particular, several results of K.~Langmann in
  \cite{Langmann:Crelle90}, \cite{Langmann:MA94}, and
  \cite{Langmann:MN00}, obtained by Diophantine approximation
  techniques, appear as special cases of some of our results.
  
  Our method is completely different. We use elementary group theory,
  valuation theory, and Siegel's theorem about integral points on
  algebraic curves. Indeed, using the Siegel-Lang extension of
  Siegel's theorem, most of our results hold over more general fields.
  \end{abstract}

\section{Introduction}
Let $k$ be a finitely generated field extension of $\QQQ$, and $R$ a
finitely generated subring of $k$. In a typical situation $k$ is a
number field, and $R=\ok$ is the ring of integers. Let $f(t,X)\in
k(t)[X]$ be an irreducible polynomial. By the well-known Hilbert
irreducibility theorem there are infinitely many specializations
$t\mapsto\bar t\in R$ such that $f(\bar t,X)$ is irreducible over
$k$. Furthermore, easy examples show that nevertheless $f(\bar t,X)$
may be reducible for infinitely many $\bar t\in R$.

Denote by $\Red_f(R)$ the set of those $\bar t\in R$ for which $f(\bar
t,X)$ is defined and reducible over $k$.

The purpose of this paper is to give several sufficient conditions
which guarantee that $\Red_f(R)$ is a finite set, and to give
non-trivial examples for infinite $\Red_f(R)$.

Our sufficient conditions are of various types. Section
\ref{SS:CondRam} gives criteria on the ramification of the place
$t\mapsto\infty$ of $k(t)$ in a root field of $f(t,X)$ which imply
finiteness of $\Red_f(R)$. For instance, if this place is not ramified
at all and $f(t,X)$ has odd degree in $X$, then $\Red_f(R)$ is
finite.

In Section \ref{SS:PolySpF} these results are applied to polynomials
of special forms. For instance, we extend a result of Langmann on Thue
polynomials. An example of this sort is the following. Let $H(t,X)\in
k[t,X]$ be a homogeneous and separable polynomial of degree $>2$. Then
$\Red_{H(t,x)-1}(R)$ is finite. This was shown by Langmann in
\cite{Langmann:MN00} under the additional assumptions that $k=\QQQ$,
$R=\ZZZ$, and that $H$ has odd degree. Our method allows to easily
obtain results about polynomial of the form $P(X)-tQ(X)$, which again
extend previous results by Langmann \cite{Langmann:Crelle90},
\cite{Langmann:MA94} by removing technical conditions he had to impose
to make his diophantine approximation techniques work.

Another sufficient condition which yields finite $\Red_f(R)$ is a
transitivity assumption on the Galois group of $f(t,X)$ over $k(t)$.
Assume that this Galois group permutes doubly transitively the roots
of $f(t,X)$. Then $\Red_f(R)$ is finite, unless $f(t,X)$ is absolutely
irreducible, and the curve $f(t,X)=0$ has genus $0$. This is shown
Section \ref{SS:DoubTr}, where we base our proof on a genus estimation
of function fields which we consider interesting in its own right.

This is the only situation where we also consider specialization in
$k$. We prove finiteness of $\Red_f(k)$ under the stronger sufficient
(and generally necessary) condition that the curve $f(t,X)=0$ has
genus $>1$.

While we obtain quite satisfactory results if the Galois group of
$f(t,X)$ is doubly transitive, the situation changes drastically if we
impose the weaker assumption that this Galois group is primitive. If a
weak condition on the composition factors is satisfied, then
$\Red_f(R)$ is finite. We also prove a converse to this criterion.

A very precise result is possible if $f(t,X)$ has prime degree in $X$,
and $k=\QQQ$, $R=\ZZZ$. If $f(t,X)=h(X)-t$ with $h(X)\in\ZZZ[X]$, then
clearly $\abs{\Red_f(\ZZZ)}=\infty$. In Section \ref{SS:PrimeDeg} we
show that this is essentially the only instance for odd prime degree
polynomials $f(t,X)$ with $\abs{\Red_f(\ZZZ)}=\infty$. It is
interesting that there are exceptions in degree $2$. This section is
related to \cite{PM:HIT}, the precursor to this paper.

Our main tool for all these results is Siegel's theorem about
algebraic curves with infinitely many integral points in a number
field, or the extension by Lang to points in a finitely generated
integral domain of characteristic $0$.

A variation of the classical reduction theorem in the proof of
Hilbert's irreducibility theorem gives the following: Assume that
$\abs{\Red_f(R)}=\infty$. Then the splitting field of $f(t,X)$ over
$k(t)$ contains an element $z$, such that $t=g(z)$ for $g(Z)\in k(Z)$
a rational function which assumes infinitely many values in $R$ on
$k$. Furthermore, $f(t,X)$ becomes reducible over $k(z)$. The property
$\abs{R\cap g(k)}=\infty$ is rather strong, results of Siegel-Lang
give precise and restrictive information about the ramification of the
places of $k(z)$ which lie above $t\mapsto\infty$.

It is clear in this setting that the Galois group of $g(Z)-t$ over
$k(t)$ is a homomorphic image of the Galois groups of $f(t,X)$. Thus
one can expect further results if one is able to classify these former
Galois groups. This has been carried out in \cite{PM:MonSiegel}. The
proof (as well as the result) is long and involved. In addition, it
makes heavy and frequent use of the classification of the finite
simple groups. The final Section \ref{S:AppClassif} provides two
applications. The first is a sufficient condition on the composition
factors of the Galois group of $f(t,X)$ if this group is primitive
which guarantees that $\Red_f(R)$ is finite.

The second result says the following: Let $f(t,X)\in\QQQ(t)[X]$ have a
simple Galois group over $\QQQ(t)$ of order $\ge3$ which is not
isomorphic to an alternating group. Then this group is preserved for
all but finitely many specializations $t\mapsto\bar t\in\ZZZ$. Even
though this is a smooth result, we doubt that a proof of that can be
achieved without knowing the list of the finite simple groups.

This work was inspired by M.\ Fried's observation of the applicability
of group theoretic methods in the analysis of Hilbert sets, see
\cite{Fried:HIT}, \cite{Fried:SC}, \cite{Fried:Sprindzuk}.

\section{Galois Theoretic Preparation}
\subsection{Description of Hilbert Sets}

The following proposition gives a convenient description of the sets
$\Red_f(R)$. The argument is a variation of the classical reduction
argument in the proof of (see e.g.\ \cite[Chapter 9]{Lang:Dio}),
combined with Lang's extension of Siegel's theorem about integral
points on algebraic curves \cite[Chapter 8]{Lang:Dio}. An alternative
argument for a similar result, which also relies on a reduction to
Siegel's Theorem, has been given by Fried, see \cite{Fried:HIT}.

\begin{Proposition}
  \label{P:Hset}
  Let $k$ be a field which is finitely generated over $\QQQ$, and $R$
  a subring of $k$ which is finitely generated over $\ZZZ$. For an
  irreducible polynomial $f(t,X)\in k(t)[X]$ of degree $\ge2$ set
\[
\Red_f(R):=\{\bar t\in R|f(\bar t,X)\text{ is defined and
  reducible}\}.
\]
Let $L$ be a splitting field of $f(t,X)$ over $k(t)$. Then there are
finitely many $z_i\in L$ and rational functions $g_i(Z)\in k(Z)$ with
$g_i(z_i)=t$, such that the following holds.
\begin{itemize}
\item[(a)] $\Red_f(R)$ and $\bigcup_i(g_i(k)\cap R)$ differ by a
  finite set.
\item[(b)] $f(t,X)$ is reducible over $k(z_i)$.
\end{itemize}
\end{Proposition}

\begin{proof}
  In order to prove the assertion, we may replace $R$ by an extension
  which still fulfills the assumption on $R$. A finitely generated
  extension of $R$ allows to assume that $k$ is the quotient field of
  $R$. Another finitely generated extension allows to assume that $R$
  is integrally closed in $k$, see \cite[Chapter 2, Prop.\ 
  4.1]{Lang:Dio}.

By replacing $X$ and $f(t,X)$ by multiples with elements in $k(t)$, we
may assume that $f(t,X)\in R[t,X]$ is monic in $X$. Let
$x_1,x_2,\dots,x_n$ be the roots of $f(t,X)$ in an algebraic closure
of $k(t)$. For each $I\subset\{1,2,\dots,n\}$ with $1\le\abs{I}\le
n-1$ set
\[
F_I(X):=\prod_{i\in I}(X-x_i).
\]
Let $K_I$ be the field generated by $k(t)$ and the coefficients of
$F_I$. Let $\beta_I$ be a primitive element of $K_I/k(t)$. We may
assume that $\beta_I$ lies in the ring generated by $R[t]$ and the
coefficients of $F_I$. In particular, $\beta_I$ is integral over
$R[t]$. Let $P_I(t,Y)\in R[t,Y]$ be the minimal polynomial of
$\beta_I$ over $k(t)$.

Now take $\bar t\in\Red_f(R)$ such that $f(\bar t,X)$ is separable.
(This assumption excludes only finitely many elements $\bar t$ from
consideration.) Write $f(\bar t,X)=u(X)v(X)$ with $u,v$ monic
polynomials in $R[X]$. As $k[t][x_1,x_2,\dots,x_n]$ is integral over
$k[t]$, the specialization map $t\mapsto\bar t$ from $k[t]$ to $k$
extends to a $k$--algebra homomorphism
$\omega:k[t][x_1,x_2,\dots,x_n]\to
k[\bar{x}_1,\bar{x}_2,\dots,\bar{x}_n]$, where the $\bar{x}_i$ are the
roots of $f(\bar t,X)$. Label these roots such that
$\omega(x_i)=\bar{x}_i$. Let $I$ be the set of those $i$ such that
$\bar{x}_i$ is a root of $u$. Denote by $\omega(F_I)$ the polynomial
$F_I$ with $\omega$ applied to its coefficients, thus $\omega(F_I)=u$.
Clearly $P_I(\bar t,\omega(\beta_I))=0$. But, by the construction
above, $\beta_I$ is a polynomial over $k[t]$ in the coefficients of
$F_I$, hence $\omega(\beta_I)\in k$, and then $\omega(\beta_I)\in R$
because $\omega(\beta_I)$ fulfills an integral equation over $R$ and
$R$ is integrally closed in $k$. Thus each such $\bar t$ gives rise to
a point $(\bar t,\omega(\beta_I))\in R^2$ on $P_I$ for some index set
$I$.

Now consider those $I$ which appear infinitely many times. Thus the
curve $P_I(T,Y)=0$ has infinitely many points with coordinates in $R$.
The Siegel--Lang Theorem \cite[Chapter 8]{Lang:Dio} implies that this
curve is rational over $k$, so there is $z_I\in L$ such that
$k(t,\beta_I)=k(z_I)$. Thus $t=g_I(z_I)$ for a rational function
$g_I\in k(Z)$. From this we obtain (a) and (b), because $F_I(X)\in
k(z_I)[X]$ is a proper factor of $f(t,X)$ over $k(z_I)$. Note that (b)
is equivalent to $f(g_I(Z),X)$ being reducible over $k(Z)$.
\end{proof}

\subsection{Not Absolutely Irreducible Polynomials}

In this section $k$ may be any field of characteristic $0$.

It is a well--known consequence from Bezout's Theorem that if
$f(X,Y)\in k[X,Y]$ is irreducible, but not absolutely irreducible,
then there are only finitely many $(a,b)\in k^2$ with $f(a,b)=0$.
Corollary \ref{C:HIT:absirr} shows that under certain additional
assumptions an analogue of this observation holds in the context of
Hilbert sets.

\begin{Lemma}\label{L:absirr}
  Let $f(t,X)\in k(t)[X]$ be an irreducible polynomial over $k$. Let
  $\ell$ be a Galois extension of $k$, and $h(t,X)\in\ell(t)[X]$ be an
  $\ell$--irreducible factor of $f(t,X)$. Then for all but finitely
  many $\bar t\in k$ with $f(\bar t,X)$ reducible over $k$, the
  polynomial $h(\bar t,X)$ is reducible over $\ell$.
\end{Lemma}

\begin{proof}
  Without loss assume that $f(t,X)\in k[t,X]$ is monic in $X$. Then we
  may assume (Gau\ss\ Lemma) that $h(t,X)$ is monic in $X$ as well.
  
  Let $\Gamma\subseteq\Gal(\ell/k)$ be chosen such that the
  $h^\gamma(t,X)$, $\gamma\in\Gamma$, are exactly the different
  conjugates of $h(t,X)$. These Galois conjugates
  $h^\gamma(t,X)\in\ell(t)[t,X]$ divide $f(t,X)$, they are relatively
  prime because they are irreducible over $\ell$, their product is
  stable under $\Gal(\ell/k)$, therefore
  $f(t,X)=\prod_{\gamma\in\Gamma}h^\gamma(t,X)$.
  
  Now let $\bar t\in k$ be such that $h(\bar t,X)$ is irreducible over
  $\ell$ and $f(\bar t,X)$ is separable. Let $g(X)\in k[X]$ be a
  non-constant factor of $f(\bar t,X)$ which is divisible by $h(\bar
  t,X)$. Then $h^\gamma(\bar t,X)$ divides $g(X)$ for each
  $\gamma\in\Gamma$. The polynomials $h^\gamma(\bar t,X)$ are
  relatively prime by separability of $f(\bar t,X)$, so $f(\bar t,X)$,
  which is the product of these polynomials, divides $g(X)$, so
  $f(\bar t,X)=g(X)$, and the claim follows.
\end{proof}

\begin{Corollary}\label{C:HIT:absirr}
  Let $f(t,X)\in k(t)[X]$ be an irreducible polynomial over $k(t)$,
  and let $A$ be the Galois group of $f(t,X)$ over $k(t)$. Assume that
  $f(\bar t,X)$ is reducible over $k$ for infinitely many $\bar t\in
  k$, and that one of the following holds.
  \begin{itemize}
  \item[(a)] $A$ is a simple group; or
  \item[(b)] $A$ acts primitively on the roots of $f(t,X)$.
  \end{itemize}
  Then $f(t,X)$ is absolutely irreducible over $k$.
\end{Corollary}

\begin{proof}
  Let $\bk$ be an algebraic closure of $k$, and $G$ the Galois group
  of $f(t,X)$ over $\bk(t)$. Suppose that $f(t,X)$ is not absolutely
  irreducible. Then $G$ is an intransitive normal subgroup of $A$.
  Hypothesis (a) implies that $G=1$. We get $G=1$ also from hypothesis
  (b), because the orbits of $G$ are a system of imprimitivity for
  $A$. Thus $f(t,X)$ decomposes into linear factors over $\bar k$,
  contrary to the previous lemma.
\end{proof}

\begin{Remark}\label{R:AbsIrr}
  The assertion of the corollary becomes false if we relax the
  assumption on $A$. For instance, take
  $f(t,X)=X^4+2(1-t)X^2+(1+t)^2$. Then $f(t,X)$ is irreducible over
  $\QQQ(t)$, but $f(t,X)=(X^2+2iX-1-t)(X^2-2iX-1-t)$, where $i^2=-1$.
  Furthermore, from $f(u^2,X)=(X^2+2uX+u^2+1)(X^2-2uX+u^2+1)$ we see
  that $f(\bar t,X)$ is reducible over $\QQQ$ for each square $\bar
  t\in\QQQ$.
\end{Remark}

\section{Siegel Functions}
\subsection{Connection with Hilbert sets}

Let $k$ be a field which is finitely generated over $\QQQ$, and $R$ a
subring which is finitely generated over $\ZZZ$. Let $f(t,X)\in
k(t)[X]$ be an irreducible polynomial over $k(t)$. Proposition
\ref{P:Hset} gives (up to finitely many exceptions) a description of
the set $\Red_f(R)$ of specializations $\bar t\in R$, such that
$f(\bar t,X)$ becomes reducible, as a union of finitely many infinite
sets of the form $g(k)\cap R$, where $g(Z)\in k(Z)$ is a rational
function. This leads to the following

\begin{Definition}
  Let $k$ be a field which is finitely generated over $\QQQ$, and
  $g(Z)\in k(Z)$ a non-constant rational function. We say that $g(Z)$
  is a \emph{Siegel function} over $k$, if there is a finitely
  generated subring $R$ of $k$ with $\abs{g(k)\cap R}=\infty$. If
  $k=\QQQ$, then we require more strongly that
  $\abs{g(\QQQ)\cap\ZZZ}=\infty$.
\end{Definition}

The condition that $g$ assumes infinitely many values in $R$ on $k$ is
quite strong, and puts severe restrictions on the form of $g$. The
basic result is due to Siegel \cite{Siegel} in the number field case,
and has been extended by Lang \cite[Theorem 8.5.1]{Lang:Dio} to the
more general fields $k$.

\begin{Proposition}\label{P:Siegel2}
  Let $k$ be a field which is finitely generated over $\QQQ$, and
  $g(Z)\in k(Z)$ a Siegel function over $k$. Then
  $\abs{g^{-1}(\infty)}\le2$.
  
  If $k=\QQQ$, and $\abs{g^{-1}(\infty)}=2$, then the two elements in
  $g^{-1}(\infty)$ are real and algebraically conjugate.
\end{Proposition}

\subsection{Cycle types of inertia generators}

Let $g(Z)\in k(Z)$ be a non--constant rational function over a field
$k$ of characteristic $0$, and $t$ a transcendental.

The following lemma is well--known (and easy to prove using Puiseux
series, for instance).

\begin{Lemma}\label{L:Puiseux}
  Let $m_1, m_2, \dots, m_r$ be the multiplicities of the elements of
  $\PP(\bk)$ in the fiber $g^{-1}(\alpha)$ for $\alpha\in\bar
  k\cup\{\infty\}$. Let $L$ be a splitting field of $g(Z)-t$ over
  $k(t)$, and $I$ the inertia group of a place of $L$ lying above the
  place $t\mapsto\alpha$ of $k(t)$. Then $I$ is cyclic, and generated
  by an element which has cycle lengths $m_1, m_2, \dots, m_r$ in the
  action on the roots of $g(Z)-t$.
\end{Lemma}

\subsection{Decomposition groups}

It is clear from Proposition \ref{P:Hset} that in order to understand
the dependency of the Hilbert sets $R\setminus\Red_f(R)$ in terms of
$A=\Gal(f(t,X)/k(t))$, one has to get control over the possibilities
for the Galois group of $g(Z)-t$ over $k(t)$ for a Siegel function
$g$.

Let $k$ be a field of characteristic $0$, and $g(Z)\in k(Z)$ be a
non-constant rational function with $\abs{g^{-1}(\infty)}\le2$. Denote
by $L$ a splitting field of $g(Z)-t$ over $k(t)$. Set
$A:=\Gal(L/k(t))$, considered as a permutation group on the roots of
$g(Z)-t$, and let $G\trianglelefteq A$ be the normal subgroup
$\Gal(\bk L/\bk(t))$.

The following lemma is a variation of the branch cycle argument, see
\cite[2.2.3]{MM} and \cite[Lemma 2.8]{H:Buch}.

\begin{Lemma}\label{L:Siegel:rat}
  Let $D\le A$ and $I\trianglelefteq D$ be the decomposition and
  inertia group of a place of $L$ lying above the place
  $t\mapsto\infty$ of $k(t)$, respectively. Then $I$ is generated by
  an element $\sigma\in I$, and the following holds.
\begin{itemize}
\item[(1)] $\sigma$ has at most two cycles, with lengths equal the
  multiplicities of the elements in $g^{-1}(\infty)$.
\item[(2)] $A=GD$ and $I\le G\cap D$.
\end{itemize}

Suppose that $k=\QQQ$, $\abs{g^{-1}(\infty)}=2$, and the two elements
in $g^{-1}(\infty)$ are real and algebraically conjugate. Then $g$ has
even degree $2m$, and the following holds.
  
\begin{itemize}
\item[(a)] $\sigma$ is a product of two $m$--cycles.
\item[(b)] $\sigma^r$ is conjugate in $D$ to $\sigma$ for all $r$ prime to
  $m$.
\item[(c)] $D$ contains an element which switches the two orbits of
  $I=\gen{\sigma}$.
\item[(d)] $D$ contains an element $\tau$ of order $1$ or $2$, such
  that $\sigma^\tau=\sigma^{-1}$, and $\tau$ fixes the orbits of $I$
  setwise.
\item[(e)] If $g^{-1}(\infty)\nsubseteq\QQQ(\zeta_m)$ (with $\zeta_m$
  a primitive $m$--th root of unity), then $D$ contains an element
  which interchanges the two orbits of $I$ and centralizes $I$.
\end{itemize}
\end{Lemma}

\begin{proof}
  Assertions (1) and (a) follow from Lemma \ref{L:Puiseux}.
  
  Assertion (2): Let $O_L$ be the valuation ring of the given place of
  $L$, and $\frakP$ the corresponding valuation ideal. Then $(O_L\cap
  k(t))/(\frakP\cap k(t))$ is naturally isomorphic to $k$. Using this
  identification, $O_L/\frakP$ is a Galois extension of $k$ with group
  $D/I$, see \cite[Chapter I, \S7, Prop.~20]{Serre:LF}. On the other
  hand, $L\cap\bk$ embeds into $O_L/\frakP$, so $D/I$ surjects
  naturally to $A/G=\Gal(L\cap\bk/k)$. Furthermore, if $\phi\in I$,
  then $u-u^\phi\in\frakP$ for all $u\in L\cap\bk$, hence $\phi$ is
  trivial on $L\cap\bk$, so $I\le D\cap G$.

It remains to prove (b) to (e).

Composing $g$ with linear fractional functions over $\QQQ$ allows to
assume that the two elements in the fiber $g^{-1}(\infty)$ are
$\pm\sqrt{d}$, where $d>1$ is a squarefree integer. Thus, without
loss, assume that $g(Z)=h(Z)/(Z^2-d)^m$, where $h(Z)\in\QQQ[Z]$ with
$\deg(h)\le 2m$, and $h(\pm\sqrt{d})\ne0$.
  
Let $y$ be a transcendental over $\QQQ$, such that $y^m=1/t$. Fix a
square root $\sqrt{d}$ of $d$, and let $\varepsilon\in\{-1,1\}$.
Substituting $y\tilde Z+\varepsilon\sqrt{d}$ for $Z$ in the equation
$h(Z)-t\cdot(Z^2-d)^m=0$ gives
\[
h(y\tilde Z+\varepsilon\sqrt{d})-\tilde Z^m(y\tilde
Z+2\varepsilon\sqrt{d})^m=0.
\]
This latter equation, by Hensel's Lemma, is solvable in the power
series ring $\bQ[[y]]$.

Thus, for $i=1,2,\dots,m$ and $\varepsilon\in\{-1,1\}$, we can
represent the $2m$ roots of $g(Z)-t$ in the form
\[
z_{i,\varepsilon}=\varepsilon\sqrt{d}+
                  a_{1,\varepsilon}\zeta^iy+
                  a_{2,\varepsilon}\zeta^{2i}y^2+
                  \dots\in\bQ[[y]],
\]
where $\zeta$ is a primitive $m$--th root of unity.

Thus $L$ can be regarded as a subfield of $\bQ((y))$. Each
automorphism of $\bQ((y))$ which fixes $y^m=1/t$ then restricts to an
element in $D\le A$, and if it is the identity on $\bQ$, then the
restriction to $L$ lies in $I$.

We will now construct suitable automorphisms of $\bQ((y))$ which, when
restricted to $L$, give the required actions on the roots of $g(Z)-t$.

To (b). Let $\hat\tau\in\Gal(\bQ((y))/\QQQ((y)))$ with
$\zeta^{\hat\tau}=\zeta^r$, and $\tau:=\hat\tau|_L$. Then
$\tau^{-1}\sigma\tau$ is the identity on $\bQ$, but
$y^{\tau^{-1}\sigma\tau}=y^{\sigma\tau}=(\zeta y)^\tau=\zeta^r y$, so
$\tau^{-1}\sigma\tau=\sigma^r$.

To (c). Choose $\hat\tau\in\Gal(\bQ((y))/\QQQ((y)))$ such that
$\sqrt{d}^{\hat\tau}=-\sqrt{d}$.

To (d). Choose $\hat\tau\in\Gal(\bQ((y))/\QQQ((y)))$, such that the
restriction of $\hat\tau$ to $\bQ$ is the complex conjugation for a
fixed embedding of $\bQ$ into $\CCC$. Then $r=-1$ in the notation of
case (b).

To (e). If $\sqrt{d}\not\in\QQQ(\zeta)$, then there is an element
$\hat\tau\in\Gal(\bQ((y))/\QQQ((y)))$ such that $\hat\tau$ moves
$\sqrt{d}$, but is the identity on $\QQQ(\zeta)$. Set
$\tau:=\hat\tau|_L$. This gives $r=1$ in case (b).
\end{proof}

\subsection{Indecomposability versus absolute indecomposability}
The main results of this paper do not depend on this section, Theorem
\ref{T:absind} below is used only in the proof of Theorem
\ref{T:mainU}.

In this section $k$ is any field of characteristic $0$. We say that a
non-constant rational function $g(Z)\in k(Z)$ is \emph{functionally
  indecomposable} if $g(Z)$ cannot be written as a composition of
rational functions in $k(Z)$ of lower degree. A classical result by
M.~Fried says that functionally indecomposable polynomials $g(Z)\in
k[Z]$ are functionally indecomposable over $\bk$,
\cite{FriedMacRae}. We extend this result to rational functions with
$\abs{g^{-1}(\infty)}\le2$.

We remark that there are many examples of functionally indecomposable
rational functions over $k$ which decompose over $\bk$. A small
example is $f(Z)=(Z^4+16Z)/(Z^3-2)$, which is indecomposable over
$\QQQ$, but decomposes over $\QQQ(\sqrt[3]{2})$. An infinite series can be
constructed a follows: Let $p$ be an odd prime, and $E$ an elliptic
curve over $\QQQ$ whose $p$-torsion points generate a field with
Galois group $\GL_2(p)$ over $\QQQ$. Denote by $[p]$ the
multiplication by $p$ map, and by $\tau$ the canonical involution on
$E$. Then $[p]$ induces a map $g:E/\gen{\tau}\to E/\gen{\tau}$. We may
interpret $g$ as a rational function, because $E/\gen{\tau}$ is a
rational curve. It is easy to see that $g$ is indecomposable over
$\QQQ$. Let $P$ be a subgroup of $E(\bk)$ of order $p$. Then $[p]$
factors as $E\to E/P\to E$, where the second isogeny is the dual of
the first one. Dividing by the canonical involution gives a
decomposition of $g(Z)$ over $\bQ$.

\begin{Theorem}\label{T:absind} Let $k$ be a field of characteristic
  $0$, and $g(Z)\in k(Z)$ be functionally indecomposable over $k$.
  Suppose that $\abs{g^{-1}(\infty)}\le2$. Then $g(Z)$ is functionally
  indecomposable over $\bk$.
\end{Theorem}

\begin{proof} The proof is by group theory. Let $L$ be a splitting
  field of $g(Z)-t$ over $k(t)$, and $\hat k$ the algebraic closure of
  $k$ in $L$. Set $A:=\Gal(L/k(t))$, and $G:=\Gal(L/\hat
  k(t))\trianglelefteq A$. Note that $G=\Gal(g(Z)-t/\bar k(t))$. The
  assumption and L\"uroth's theorem give that $A$ is primitive on the
  roots of $g(Z)-t$. So it remains to show that $G$ is primitive as
  well. Let $I$ and $D$ be the inertia and decomposition group of a
  place of $L$ lying above $t\mapsto\infty$.  Lemma \ref{L:Siegel:rat}
  gives $A=GD$, and $I$ has at most two orbits.  Thus the theorem
  follows from the following purely group theoretic result.\end{proof}

\begin{Theorem}\label{T:agdi} Let $\Omega$ be a finite set, and let
  $A\le\sym(\Omega)$ be a primitive permutation group on $\Omega$. Let
  $1<G\trianglelefteq A$ be a normal subgroup, which contains a cyclic
  subgroup $I$ with the following properties:\begin{itemize}
\item[(a)] $I$ has at most two orbits on $\Omega$, and
\item[(b)] $A=GN_A(I)$, where $N_A(I)$ denotes the normalizer of $I$
  in $A$.
\end{itemize}
Then $G$ acts primitively on $\Omega$ as well.
\end{Theorem}

\begin{proof} Suppose that $G$ acts imprimitively. Then $\Omega$
  is a disjoint union of $\Delta=\Delta_1$, $\Delta_2$, \ldots,
  $\Delta_m$, where $1<r=\abs{\Delta_i}<\abs{\Omega}$, and $G$
  permutes the $m$ sets $\Delta_i$. We assume that among these systems
  we have chosen one such that $\abs{\Delta}$ is maximal. This implies
  that $G$ permutes the $\Delta_i$'s primitively.
  
  If $a\in A$, then the sets $\Delta_i^a$, $i=1,\dots,m$, again
  constitute a system of imprimitivity for $G$, this follows from
  $G\trianglelefteq A$. We claim that there is an element $a\in A$
  such that $\Delta^a$ is not contained in an $I$-orbit. Suppose that
  is not the case. Then, for each $a$, each orbit of $I$ is a union of
  sets $\Delta_i^a$. So the sets $\Delta_i^a$ in an orbit of $I$ are
  the orbits of a subgroup of $I$. The size of this subgroup in $I$
  depends on $\abs{\Delta}$, but not on $a$. On the other hand, a
  subgroup in a cyclic group is uniquely given by its order. We obtain
  that for each $a\in A$ the sets $\Delta_i^a$ are a permutation of
  the sets $\Delta_i$, thus the $\Delta_i$ are a system of
  imprimitivity for $A$, contrary to the assumption that $A$ is
  primitive.
  
  Thus $I$ has two orbits, and we may assume that $\Delta$ intersects
  them both non-trivially. So $I$ permutes the $\Delta_i$
  transitively. Let $K\triangleleft G$ be the kernel of the action of
  $G$ on the $\Delta_i$, and let $I_\Delta$ the setwise stabilizer of
  $\Delta$ in $I$. As $IK/K$ permutes the $\Delta_i$ regularly, we
  obtain that $I_\Delta$ fixes each $\Delta_i$ setwise, so
  $I_\Delta\le K$. As $A$ is primitive, and $K$ is intransitive on
  $\Omega$, we get $\bigcap_{a\in A}K^a=1$. From $A=GN_A(I)$ and
  $K\triangleleft G$ we obtain $\bigcap_{a\in N_A(I)}K^a=1$. But
  $I_\Delta\le K$ and $(I^\Delta)^a=I_\Delta$ for all $a\in N_A(I)$,
  so $I_\Delta=1$. On the other hand, $I_\Delta$ has two orbits on
  $\Delta$, so this implies $\abs{\Delta}=2$.
  
  Choose $\delta\in\Delta$, and let $A_\delta$ and $G_\delta$ be the
  stabilizers of $\delta$ in $A$ and $G$, respectively. Also, let
  $G_\Delta$ be the setwise stabilizer of $\Delta$. Clearly
  $[G_\Delta:G_\delta]=2$, so $G_\delta$ is normal in $G_\Delta$.
  Furthermore, $G_\delta=A_\delta\cap G$ is normal in $A_\delta$, so
  $G_\delta$ is normal in the group $U:=\gen{A_\delta,G_\Delta}$. But
  $A_\delta$ is a maximal subgroup of $A$ by primitivity of $A$, so
  $U=A_\delta$ or $U=A$. The former possibility cannot hold, because
  $G_\Delta$ is transitive on $\Delta$, so $G_\Delta\le A_\delta$
  cannot hold. Thus $U=A$, so $G_\delta$ is normal in $A$, hence
  $G_\delta=1$. So $G$ acts regularly on $\Omega$, $G_\Delta=K$, and
  $G$ is the direct product of $G_\Delta$ and $I$ by order reasons.
  But then the intransitive group $I$ is normal in $A=GN_A(I)$,
  contrary to primitivity of $A$.
\end{proof}

\begin{Remark} If $I$ has only one orbit on $\Omega$, then we got the
  claim without using assumption (b).
  
  However, in general we cannot remove the assumption (b), there are
  infinite series of counterexamples. For instance let $m\ge3$ be an
  integer, and $A=(\sym_m\times\sym_m)\rtimes C_2$, where $C_2$ flips
  the two components. Let the action be given on the coset space
  $A/A_1$, where $A_1=(\sym_{m-1}\times\sym_{m-1})\rtimes C_2$. This
  action is easily be seen to be primitive. Let
  $G=\sym_m\times\sym_m$, and $I$ be generated by $(a,b)$, where $a$
  is an $m$-cycle, and $b$ is an $(m-1)$-cycle. Then one verifies that
  $I$ has two orbits. However, $G$ is not primitive anymore, because
  $\sym_{m-1}\times\sym_m$ is a group properly between $G_1=G\cap A$
  and $G$.
\end{Remark}

\section{Consequences from $\abs{\Red_f(R)}=\infty$.}

Let $k$ be a field which is finitely generated over $\QQQ$, and $R$ a
finitely generated subring. The following lemma summarizes how we use
the information that $\abs{\Red_f(R)}$ is an infinite set for an
irreducible polynomial $f(t,X)\in k(t)[X]$.

\begin{Lemma}\label{L:RedGal} Let $f(t,X)$ be irreducible, and assume
  that $\Red_f(R)=\infty$. Let $L$ be a splitting field of $f(t,X)$
  over $k(t)$, and $x\in L$ a root of $f(t,X)$. Set $A:=\Gal(L/k(t))$,
  and let $D$ and $I$ be the decomposition group and inertia group of
  a place of $L$ lying above $t\mapsto\infty$, respectively.
  
  Then there is a rational field $k(t)\subseteq k(z)\subseteq L$, such
  that the following holds, where $A_x$ and $A_z$ are the stabilizers
  in $A$ of $x$ and $z$, respectively.

\begin{enumerate}
\item[(a)] $A_z$ acts intransitively on the coset space $A/A_x$.
\item[(b)] $I$ is cyclic, and has at most two orbits on $A/A_z$. If
  $k=\QQQ$ and $R=\ZZZ$, then these orbits have equal lengths.
\item[(c)] If $k=\QQQ$ and $R=\ZZZ$, then $D$ is transitive on
  $A/A_z$.
\end{enumerate}
\end{Lemma}

\begin{proof} The existence of the field $k(z)$ with (a) follows from
  Proposition \ref{P:Hset}. Furthermore, $t=g(z)$ with $g(Z)\in k(Z)$ 
  a Siegel function, so (b) and (c) follow from Lemma
  \ref{L:Siegel:rat}(1), (a), and (c).
\end{proof}

\section{Applications}

There are several results by K.~Langmann \cite{Langmann:Crelle90},
\cite{Langmann:MA94}, \cite{Langmann:MN00}, where he studies integral
Hilbert sets of irreducible polynomials $f(t,X)\in\QQQ[t,X]$, if $f$
assumes a very specific form. The types he considers are so-called
Thue-equations, where $f(t,X)=H(t,X)-1$ with $H$ a homogeneous
polynomial, or polynomials of the form $P(X)-tQ(X)$, or modifications.
Under further technical assumptions, he shows that $f(t,X)$ stays
irreducible for almost all integral specializations of $t$.

Below we give irreducibility theorems which immediately imply several
generalizations of Langmann's results. Our proofs are completely
different from his. We believe that the Galois theoretic preparation
from the previous sections is the suitable setting for such
irreducibility results.

In contrast to other results given in later sections, we need very
little group theoretic techniques here.

Throughout this section $k$ denotes a field which is finitely
generated over $\QQQ$, and $R$ is a subring of $k$ which is finitely
generated over $\ZZZ$.

\subsection{Conditions on ramification}\label{SS:CondRam}
In this section we obtain finiteness results under suitable conditions
on the ramification indices of the places of a root field of $f(t,X)$
which lie above the place $t\mapsto\infty$ of $k(t)$. A proof will be
given later.

\begin{Theorem}\label{T:unram}
Let $f(t,X)\in k(t)[X]$ be irreducible, and assume that the place
$t\mapsto\infty$ of $k(t)$ is unramified in the field $k(t,x)$, where
$x$ is a root of $f$. Then one of the following holds.
\begin{itemize}
\item[(i)] $f(\bar t,X)$ is irreducible over $k$ for all but finitely
  many $\bar t\in R$, or
\item[(ii)] There is an element $z\in k(t,x)$, such that $t=g(z)$ with 
  $g(Z)\in k(Z)$ of degree $2$.
\end{itemize}
\end{Theorem}

\begin{Remark*} In general one cannot avoid the situation of case
  (ii). For instance set $g(Z)=1/(Z^2-d)$, where $d>1$ is a squarefree
  integer. Let $z$ be a root of $g(Z)-t$, and let $x$ be algebraic
  over $k(z)$ such that the places $z\mapsto\pm\sqrt{d}$ of $k(z)$ are
  unramified in $k(z,x)$. Then the minimal polynomial $f(t,X)$ of $z$
  over $k(t)$ fulfills the assumptions of the theorem.  However, there
  are infinitely many $\bar z\in\QQQ$ with $\bar t=g(\bar z)\in\ZZZ$,
  and for each such $\bar t$ the polynomial $f(\bar t,X)$ is
  reducible.
\end{Remark*}

Assume the situation of the previous theorem, and let $\tilde x$ be a
primitive element of the normal closure of $k(t,x)/k(t)$. Apply the
theorem to the minimal polynomial of $\tilde x$ over $k(t)$. (Note
that $t\mapsto\infty$ is unramified in this normal closure too.) Then
case (ii) can only appear if the Galois group of $k(t,\tilde x)/k(t)$
has a subgroup of index $2$. Thus we obtain the following

\begin{Corollary}\label{C:unram}
Let $f(t,X)\in k(t)[X]$ be irreducible, and assume that the place
$t\mapsto\infty$ of $k(t)$ is unramified in the field $k(t,x)$, where
$x$ is a root of $f$. Suppose that the Galois group $A$ of $f(t,X)$
over $k(t)$ has no subgroup of index $2$. Then $A=\Gal(f(\bar t,X)/k)$ 
for all but finitely many $\bar t\in R$.
\end{Corollary}

In general one cannot relax the assumption about the infinite place
without introducing severe other conditions in the theorem. However,
if the base field is $k=\QQQ$, then the following holds.

\begin{Theorem}\label{T:gcd=1}
  Let $f(t,X)\in\QQQ(t)[X]$ be irreducible of odd degree, and $x$ a
  root of $f$.  Assume that the greatest common divisor of the
  ramification indices of the places of $\QQQ(t,x)$ which lie above
  the place $t\mapsto\infty$ of $\QQQ(t)$ is $1$. Then $f(\bar t,X)$
  is irreducible over $\QQQ$ for all but finitely many $\bar
  t\in\ZZZ$.
\end{Theorem}

A theorem of a similar flavor is

\begin{Theorem}\label{T:rat}
  Let $f(t,X)\in\QQQ[X]$ be irreducible, and $x$ a root of $f$. Assume
  that $\QQQ(t,x)$ has a rational unramified place above the place
  $t\mapsto\infty$ of $\QQQ(t)$. Then $f(\bar t,X)$ is irreducible
  over $\QQQ$ for all but finitely many $\bar t\in\ZZZ$.
\end{Theorem}

\begin{proof}[Proof of the theorems.]
  Let us assume that $f(\bar t,X)$ is reducible for infinitely many
  $\bar t\in R$. Let $L$ be a splitting field of $f(t,X)$ over $k(t)$,
  and $x\in L$ a root of $f(t,X)$. We make frequent use of Lemma
  \ref{L:RedGal} and the notation introduced there.
  
  Let $\sigma\in A$ be a generator of $I$.
  
  First assume the situation from Theorem \ref{T:unram}. This means
  that the inertia group $I$ is trivial, so $\sigma=1$. On the other
  hand, $\sigma$ has at most two cycles on $A/A_z$. As $A_z$ is a
  proper subgroup of $A$ (because $A_z$ is intransitive on $A/A_x$),
  this implies $[A:A_z]=2$. Furthermore, $A_z$ is normal in $A$, so
  $A_zA_x$ is a proper subgroup of $A$. This implies $A_x\subseteq
  A_z$, so $z\in k(t,x)$, and the claim follows.
  
  Next assume the situation of Theorem \ref{T:rat}. The decomposition
  group $D$ acts transitively on $A/A_z$. On the other hand, the
  assumption of the rational unramified place of $\QQQ(t,x)$ above
  $t\mapsto\infty$ implies that $D$ has a fixed point on $A/A_x$. Thus
  $D$ is a subgroup of a conjugate of $A_x$. But $D$ is transitive on
  $A/A_z$, so also $A_x$ is transitive on $A/A_z$, a contradiction.
  
  Finally assume the assumptions from Theorem \ref{T:gcd=1}. Then
  $\sigma$ acts on $A/A_z$ as a product of $r$ cycles of length $m$,
  with $r=1$ or $2$.  For each place $\frakP_i$ of $\QQQ(t,x)$ above
  $t\mapsto\infty$, let $e_i$ and $f_i$ be the ramification index and
  residue degree of $\frakP_i$, respectively. It follows that $\sigma$
  has $f_1$ cycles of length $e_1$, $f_2$ cycles of length $e_2$,
  \dots, on $A/A_x$. Thus the greatest common divisor of the cycle
  lengths $a_1,a_2,\dots,a_j$ of $\sigma$ on $A/A_x$ is $1$.
  
  Recall that $A_x$ acts intransitively on $A/A_z$. Let $u<rm=[A:A_z]$
  be an orbit length of this action. As $\sigma^{a_i}$, $1\le i\le j$
  has a fixed point on $A/A_x$, it is conjugate to an element in
  $A_x$.  On the other hand, $\sigma^{a_i}$ has cycle lengths
  $m/\gcd(m,a_i)$ on $A/A_z$.  Therefore $m/\gcd(m,a_i)$ divides $u$
  and hence $\gcd(u,m)$ too. Thus $m/\gcd(m,u)$ divides $\gcd(m,a_i)$
  for each $a_i$. But the $a_i$ have the greatest common divisor $1$,
  hence $m=\gcd(m,u)$, so $r=2$ and $u=m$. In particular, $A_x$ has
  two orbits of equal length $m$ on $A/A_z$. The orbit lengths of
  $A_x$ on $A/A_z$ are proportional to the sizes of the double cosets
  $A_zaA_x$, $a\in A$. Therefore $A_z$ has two orbits of equal length
  on $A/A_x$, so $[A:A_x]=\deg_X(f(t,X))$ is even, which proves the
  theorem.
\end{proof}

\subsection{Polynomials of special forms}\label{SS:PolySpF}

The following theorem is a generalization of \cite[Satz
3.5]{Langmann:MN00}. Langmann obtains his result under the following
three additional assumptions none of which we need in our
approach:\begin{itemize}
\item[(a)] $k=\QQQ$ and $R=\ZZZ$,
\item[(b)] the degree of $H$ is odd, and
\item[(c)] $t$ does not divide $H(t,X)$.
\end{itemize}

\begin{Theorem}\label{T:Langmannsep} Let $H(t,X)\in k[t,X]$ be a
  homogeneous polynomial of total degree $>2$ which is separable with
  respect to $X$. Then $H(\bar t,X)-1$ is irreducible for all but
  finitely many $\bar t\in R$.
\end{Theorem}

\begin{Remark*} The theorem is false for degree $2$, even over the
  rationals. To see this set $H(t,X)=X^2-dt^2$ with $d>1$ a
  square--free integer. This is indeed an exception, because the
  Pellian equation $X^2-dt^2=1$ has infinitely many integral
  solutions.
\end{Remark*}

A different generalization of Langmann's result is obtained simply by
removing the separability assumption on $H(t,X)$ and replacing it by
the obviously necessary condition that $H(t,X)$ is not a proper power.

\begin{Theorem}\label{T:Langmannirr} Let $H(t,X)\in\QQQ[t,X]$ be a
  homogeneous polynomial of odd degree which is not divisible by $t$.
  If $H(t,X)$ is not a proper power in $\QQQ[t,X]$, then $H(\bar
  t,X)-1$ is irreducible for all but finitely many $\bar t\in\ZZZ$.
\end{Theorem}

\begin{Remark*} This theorem is no longer true for number fields. An
  example is the following: Let $k$ be a number field with an infinite
  group of units, and $R$ the ring of integers.  Set
  $H(t,X):=X^2(X-t)$. From
\[
H(\frac{1-Z^3}{Z},X)-1 = (X-\frac{1}{Z})(X^2+Z^2X+Z)
\]
and the fact that $\bar t=(1-\bar z^3)/\bar z\in R$ for each unit $\bar
z$ we obtain reducibility of $H(\bar t,X)-1$ for infinitely many $\bar
t\in R$.
\end{Remark*}

\begin{Conjecture} The assumption that $H(t,X)$ has odd degree $n$
  in the above theorem can be dropped if we require the following
  necessary conditions:
\begin{itemize}
\item[(a)] $n\ne2,4$.
\item[(b)] If $4$ divides $n$, then $-4H(t,x)$ is not a $4$-th power
  in $\QQQ[t,X]$. (For otherwise $H(t,X)-1$ is already reducible).
\end{itemize}
\end{Conjecture}

\begin{Remark*} The group theory got quite involved in an attempt to
  prove this conjecture. While we feel that we got close to a proof,
  some difficulties could not be settled. The conjecture is true up to
  degree $25$, at least under the slightly stronger condition that
  $H(t,X)$ is not a power in $\bar\QQQ[t,X]$. From above we know
  already that we have to assume $n\ne2$. The following example shows
  that $n\ne4$ is also a necessary condition. This is interesting
  because the associated curve has genus $1$, so the polynomial has a
  linear factor for only finitely many integral specializations:

Let $d>1$ be a squarefree integer, and set
\[
f(t,X)=-4dX^2(dX^2-t^2)-1.
\]
Note that
\[
f(\frac{Z^2+d}{Z^2-d},X)=
-(2dX^2-\frac{4dZ}{Z^2-d}X-1)(2dX^2+\frac{4dZ}{Z^2-d}X-1).
\]
There are infinitely many integers $u,v$ with $u^2-dv^2=1$. For $\bar
z=u/v$ we obtain $\bar t=\frac{\bar z^2+d}{\bar
  z^2-d}=u^2+dv^2\in\ZZZ$, and by the above factorization $f(\bar
t,X)$ is reducible.
\end{Remark*}

The proof of Theorem \ref{T:Langmannsep} is based on

\begin{Proposition}\label{P:t^mh} Let $m$ be a positive integer, and
  $h(X)\in k[X]$ a non-constant separable polynomial. Suppose that
  ${\bar t}^mh(X)-1$ is reducible for infinitely many $\bar t\in R$.
  Then $m\le2$ and $\deg(h)$ is even.
\end{Proposition}

\begin{proof}
  $t^mh(X)-1$ is irreducible over $\bk$, for instance by the
  Eisenstein criterion with respect to a linear factor of $h(X)$. Let
  $x$ be a root of $t^mh(X)-1$. By the separability of $h$ and
  Hensel's Lemma, we can write $x$ as a Laurent series in $1/t^m$ over
  $\bk$. Thus the place $t\mapsto\infty$ is unramified in $k(t,x)$,
  and so is the place $t^m\mapsto\infty$ of $k(t^m)$ in $k(x)$.
  Theorem \ref{T:unram} gives $z\in k(t,x)$ such that $k(z)$ is a
  quadratic extension of $k(t)$, in particular, $h(t,X)$ has even
  degree in $X$.
  
  Let $k(y)$ be the intersection of $k(z)$ and the normal closure of
  $k(x)/k(t^m)$.  Clearly, $t^mh(X)-1$ is reducible over $k(y)$, so
  $k(y)$ is a proper extension of $k(t^m)$.
  
  From now on we consider only the fields between $k(z)$ and $k(t^m)$.
  To ease language, we extend the coefficients to $\bk$. The place
  $t^m\mapsto\infty$ is totally ramified in $\bk(t)$, so there are at
  most two places of $\bk(y)$ above $t^m\mapsto\infty$. On the other
  hand, the place $t^m\mapsto\infty$ is unramified in the normal
  closure of $k(x)/k(t^m)$, so it is unramified in $\bk(y)$ as well.
  
  Thus $[\bk(y):\bk(t^m)]=2$, there are two places of $\bk(y)$ above
  $t^m\mapsto\infty$, and these two places are the only places which
  are ramified in $\bk(z)$. Let $\frakp$ be a place of $\bk(y)$ lying
  above $t^m\mapsto0$. From what we saw, $m$ places of $\bk(z)$ lie
  above $\frakp$. Thus at least $m$ places of $\bk(z)$ lie above
  $t^m\mapsto0$. On the other hand, $t^m\mapsto0$ is totally ramified
  in $\bk(t)$, so at most two places of $\bk(z)$ lie above
  $t^m\mapsto0$. Thus $m\le2$, and the claim follows.
\end{proof}

\begin{proof}[Proof of Theorem \ref{T:Langmannsep}]
  Let $n$ be the total degree of $H(t,X)$. Then $H(t,X)=t^nh(X/t)$,
  where $h(X)\in k[X]$ is a polynomial of degree $\le n$. Note that
  $H(\bar t,X)-1$ is reducible if and only if $\bar t^nh(X)-1$ is
  reducible, so the claim follows from Proposition \ref{P:t^mh}.
\end{proof}

\begin{proof}[Proof of Theorem \ref{T:Langmannirr}] Let $n$ be the
  total degree of $H(t,X)$, and $e$ the greatest common divisor of the
  multiplicities of the linear factors of $H(t,X)$ over $\bar\QQQ$.
  Then $H(t,X)=\tilde H(t,X)^e$, where $\tilde H(t,X)\in\bar\QQQ[t,X]$
  is homogeneous of degree $n/e$, and the greatest common divisor of
  the multiplicities of the linear factors of $\tilde H(t,X)$ is $1$.
  By Capelli's Theorem, $\tilde H(t,X)-1$ is irreducible over
  $\bar\QQQ$.
  
  Suppose that $H(\bar t,X)-1$ is reducible for infinitely many $\bar
  t\in\ZZZ$. Then there is a Siegel function $g(Z)\in\QQQ(Z)$ such
  that $H(g(Z),X)-1$ is reducible over $\QQQ(Z)$. Let $A(Z,X)$ be a
  non-trivial factor. We claim that $\tilde H(g(Z),X)-1$ is reducible
  over $\bar\QQQ(Z)$. Suppose that is not the case. As $\tilde
  H(t,X)-1$ divides $H(t,X)-1$, we may assume that $\tilde
  H(g(Z),X)-1$ divides $A(Z,X)$. However, the Galois group
  $\Gal(\bar\QQQ/\QQQ)$ fixes $A(Z,X)$, while it permutes transitively
  the factors $\tilde H(g(Z),X)-\zeta^i$ (up to scalar multiples) of
  $H(g(Z),X)-1$, where $\zeta$ is a primitive $e$-th root of unity,
  and $i=1,2,\ldots,e$. Thus $H(g(Z),X)-1$ divides $A(Z,X)$, a
  contradiction.
  
  So $\tilde H(g(Z),X)-1$ is reducible over $\bar\QQQ(Z)$ with $g(Z)$
  a Siegel function over $\QQQ$. Set $\tilde n=n/e=\deg(\tilde H)$,
  and write $\tilde H(t,X)=t^{\tilde n}h(X/t)$. Upon replacing $X$ by
  $Xg(Z)$, we get that $g(Z)^{\tilde n}h(X)-1$ is reducible over
  $\bar\QQQ(Z)$. Let $L$ be a splitting field of $th(X)-1$ over
  $\bQ(t)$, and $z$ be a root of $g(Z)^{\tilde n}-t$. So $th(X)-1$ is
  reducible over $\bQ(z)$. Denote by $\bQ(y)$ the intersection of
  $\bQ(z)$ with $L$. Of course, $th(X)-1$ is reducible over $\bQ(y)$
  as well. Write $t=\tilde g(y)$ with $\tilde g(Y)\in\bQ(Y)$. As
  $\tilde g$, composed with another rational function, gives $g$, we
  obtain that the fiber $\tilde g^{-1}(\infty)$ contains at most two
  elements, and that the multiplicities of these elements are the
  same.
  
  By construction $h(X)$ is a polynomial where the multiplicities of
  the roots have no common divisor $>1$. These multiplicities are
  exactly the ramification indices of the places of $\bQ(t,x)$ which
  lie above the place $t\mapsto\infty$ of $\bar\QQQ(t)$, where $x$ is
  a root of $th(X)-1$. Thus the assumptions of Theorem \ref{T:gcd=1}
  are fulfilled except that we are not necessarily over the rationals.
  Nevertheless, the proof of that theorem covers our situation,
  because we used the assumption that the base field is $\QQQ$ only to
  guarantee that, in the present context, the elements in the fiber
  $\tilde g^{-1}(\infty)$ have the same multiplicities. But we have
  verified this property above, so the claim follows.
\end{proof}

Another easy consequence of Theorem \ref{T:unram} (and its proof) is

\begin{Theorem}\label{T:LangmGen}
  Let $P(X)\in k[X]$ be a polynomial which is relatively prime to the
  separable polynomial $Q(X)\in k[X]$ of degree $\ge\deg(P)-1$. Then
  one of the following holds:
\begin{itemize}
\item[(a)] $P(X)-\bar tQ(X)$ is irreducible for all but finitely many
  $\bar t\in R$, or
\item[(b)] $\text{max}(\deg(P),\deg(Q))$ is even, and there is a
  rational function $g(Z)\in\bk[Z,\frac{1}{Z}]$ of degree $2$, such
  that $P(X)-tQ(X)$ factors over $\bk(Z)$ in two factors of equal
  degree in $X$.
\end{itemize}
\end{Theorem}

\begin{Remark*} This result generalizes
  \cite[Folgerung 6]{Langmann:Crelle90}, where this is proven under
  the assumption that $\deg(Q)=\deg(P)-1$. Also, the rather technical
  result \cite[Folgerung 3.4]{Langmann:MA94} is a very special case of
  Theorem \ref{T:unram}.
\end{Remark*}

A direct application of Theorem \ref{T:rat} to polynomials of the form
$P(X)-tQ(X)$ (which are studied in \cite{Langmann:Crelle90} and
\cite{Langmann:MN00}, too) is

\begin{Theorem}\label{T:PQ} Let $P(X),Q(X)\in\QQQ[X]$ be relatively prime
  polynomials, and assume that $Q(X)$ has a simple rational root. Then
  $P(X)-\bar tQ(X)$ is irreducible for all but finitely many $\bar
  t\in\ZZZ$.
\end{Theorem}

\begin{Remark*} Langmann and other authors, in particularly D{\`e}bes
  (see \cite{Debes:um}) and Fried (see \cite{Fried:Sprindzuk}) have
  studied irreducibility questions when specializing $t$ in certain
  subsets of the integers. Examples are the sets of prime powers, or
  powers of a fixed integer. A recent result of this kind with a
  completely elementary and elegant proof (in particularly not relying
  on Siegel's Theorem) is the following by Cavachi \cite{Cavachi} (his
  version is slightly more general): Let $P(X),Q(X)\in\QQQ[X]$ be
  relatively prime polynomial with $\deg(P)<\deg(Q)$. Then
  $P(X)-pQ(X)$ is irreducible for all but finitely many prime numbers
  $p$.
\end{Remark*}

\subsection{Doubly transitive Galois groups}\label{SS:DoubTr}

The main result of this section is

\begin{Theorem}\label{T:doubly} Let $f(t,X)\in k(t)[X]$ be
  irreducible, and assume that the Galois group of $f(t,X)$ over
  $k(t)$ acts doubly transitively on the roots of $f$. If $f(\bar
  t,X)$ is reducible for infinitely many $\bar t\in R$, then the
  following holds, where $x$ is a root of $f(t,X)$:
\begin{enumerate}
\item[(a)] $f(t,X)$ is absolutely irreducible, and
\item[(b)] $k(t,x)$ has genus $0$, and
\item[(c)] there are at most two places of $\bk(t,x)$ above
  $t\mapsto\infty$.
\end{enumerate}
\end{Theorem}

As a preparation we need a bound on the genus of function fields.

\subsubsection{Genus comparison}

Let $k$ be a field of characteristic $0$, and $L/k(t)$ be a finite
Galois extension. Let $\hk$ be the algebraic closure of $k$ in $L$.
Set $A:=\Gal(L/k(t))$ and $G:=\Gal(L/\hk(t))\trianglelefteq A$.

The following is well known (see e.g.\ \cite[Exp.\ XIII, Cor.\ 
2.12]{SGA1}): Let $\frakp_i$, $i=1,\ldots,r$, be the places of
$\bk(t)$ which ramify in $\bk L$. Let $I_i$ be the inertia group of a
place of $\bk L$ lying above $\frakp_i$. We identify $\Gal(\bk
L/\bk(t))$ with $G$ via restriction to $L$. We can choose elements
$\sigma_i\in G$ such that each $\sigma_i$ is conjugate to a generator
of $I_i$, and the following holds:
\begin{itemize}
\item[(a)] The $\sigma_i$, $i=1,\ldots,r$, generate $G$.
\item[(b)] $\sigma_1\sigma_2\cdots\sigma_r=1$.
\end{itemize}
The (not uniquely given) tuple $(\sigma_1,\sigma_2,\dots,\sigma_r)$ is 
called a \emph{branch cycle description} in $G$.

If $E$ is a field between $k(t)$ and $L$ with $n=[E:k(t)]$, then $A$
acts as a permutation group on the $n$ conjugates of a primitive
element of $E/k(t)$. Let $\pi_E$ be the homomorphism from $A$ to the
transitive permutation group of degree $n$.

For a permutation $\sigma$ on $n$ letters let $\ind(\sigma)$ be ``$n$
minus the number of cycles'' of $\sigma$. Let $E$ be as above, and
assume that $\hk\cap E=k$.

The Riemann-Hurwitz genus formula allows to compute the genus $g(E)$
of $E$:

\begin{equation}\label{E:RH}
2(n-1+g(E))=\sum_{i=1}^r\ind(\pi_E(\sigma_i))
\end{equation}

Associated to $\pi_E$ is the permutation character $\chi_E$, where
$\chi_E(\sigma)$ is the number of fixed points of $\pi_E(\sigma)$. In
the following lemma a character is understood as a character over the
complex numbers.

\begin{Lemma}\label{L:EF}
In the setting from above, let $F$ be another field between $k(t)$ and 
$L$, such that $\hk\cap F=k$. Suppose that $\pi_F-\pi_E$ is a
character. Then the following holds.
\begin{itemize}
\item[(a)] $g(E)\le g(F)$.
\item[(b)] For each subgroup $U\le A$, the number of orbits of
  $\pi_E(U)$ is not bigger than the number of orbits of $\pi_F(U)$.
\end{itemize}
\end{Lemma}

\begin{Remark*} To my knowledge part (a) has first been observed by
  R.~Guralnick some years ago. His proof in \cite{Gural:MG} does not
  rely on the Riemann-Hurwitz formula and the branch cycle
  description.  Instead, he uses Jacobians of function fields and the
  action of the Galois group on the $\ell$-torsion points for a
  suitable $\ell$.  This approach proves (a) in positive
  characteristic as well. Independently I had found this result by
  using a linear algebra result of Scott (see below). As Guralnick's
  proof is not yet published, we supply our elementary proof. This
  proof, however, does not work in positive characteristic due to the
  lack of branch cycle descriptions.
\end{Remark*}

The following proposition is an immediate consequence of Scott's
result \cite[Theorem 1]{Scott} and Maschke's theorem.

\begin{Proposition}
  Let the finite group $G$ act linearly on the $n$--dimensional
  complex vector space $V$. For $M$ an element or subgroup of $G$, let 
  $d(M)$ be the dimension of the subspace of fixed vectors under
  $M$. Let $G$ be generated by $\sigma_1,\sigma_2,\ldots,\sigma_r$,
  and assume that $\sigma_1\sigma_2\cdots\sigma_r=1$. Then
\[
2(n-d(G))\le\sum_{i=1}^r(n-d(\sigma_i)).
\]
\end{Proposition}

\begin{proof}[Proof of Lemma \ref{L:EF}] Let $V_E$ and $V_F$ be the
  permutation modules corresponding to $\pi_E$ and
  $\pi_F$. Considering the set which $\pi_E(G)$ acts on as the natural
  basis of $V_E$, we may consider $\pi_E$ as a homomorphism from $G$
  to $\GL(V_E)$. With respect to this natural basis, we see the
  following: If $\pi_E(\sigma)$ has a cycle of length $m$, then the
  eigenvalues of $\pi_E(\sigma)$ on the space spanned by these $m$
  cyclically moved elements are just the $m$-th roots of unity. In
  particular, the eigenvalue $1$ appears exactly once on this
  subspace. Thus the number of cycles of $\pi_E(\sigma)$ equals
  $d(\pi_E(\sigma))$. From that we obtain
\begin{equation}\label{E:E}
2([E:k(t)]-1+g(E))=\sum_{i=1}^r([E:k(t)]-d(\pi_E(\sigma_i)))
\end{equation}
and likewise
\begin{equation}\label{E:F}
2([F:k(t)]-1+g(F))=\sum_{i=1}^r([F:k(t)]-d(\pi_F(\sigma_i))).
\end{equation}
The assumption that $\pi_F-\pi_E$ is a character implies that $V_F$
has a $G$-submodule which is $G$-isomorphic to $V_E$. By Maschke's
theorem, there is a $G$--invariant complement $W$. Let
$\pi:G\to\GL(W)$ be the associated homomorphism. As $\pi_E$ and
$\pi_F$ are transitive, they both contain the principal character
$1_G$ with multiplicity $1$. Therefore $d(\pi(W))=0$. Note that
$\dim(W)=[F:k(t)]-[E:k(t)]$. The proposition gives
\begin{equation}\label{E:EF}
2([F:k(t)]-[E:k(t)])\le\sum_{i=1}^r([F:k(t)]-[E:k(t)]-d(\pi(\sigma_i))).
\end{equation}
Clearly $d(\pi_F(\sigma))-d(\pi_E(\sigma))=d(\pi(\sigma))$, so (a)
follows from \eqref{E:E}, \eqref{E:F}, and \eqref{E:EF}.

Claim (b) is obvious, because the number of orbits of $\pi_E(U)$ is the
multiplicity of the principal character $1_U$ in the restriction of
$\pi_E$ to $U$.
\end{proof}

\begin{proof}[Proof of Theorem \ref{T:doubly}.] A doubly transitive
  permutation group is primitive, so (a) follows from Corollary
  \ref{C:HIT:absirr}.
  
  Again, choose $z\in L$, where $L$ is a splitting field of $f(t,X)$
  over $k(t)$, such that $t$ is a Siegel function in $z$, and $A_x$ is
  intransitive on $A/A_z$. Let $\pi_x$ and $\pi_z$ be the permutation
  characters of the action of $A$ on $A/A_x$ and $A/A_z$,
  respectively. The scalar product $(\pi_x,\pi_z)$ of characters is
  the number of orbits of $A_x$ on $A/A_z$, so $(\pi_x,\pi_z)\ge2$.
  Each of these characters contains the principal character $1_A$ with
  multiplicity $1$. Furthermore, $\pi_x-1_G$ is irreducible, because
  $A$ is doubly transitive on $A/A_x$ (see \cite[Chapter 4, Theorem
  3.4]{Gorenstein}). Thus, as the irreducible characters are an
  orthonormal basis of the class functions on $A$, we obtain that the
  nonprincipal part of $\pi_x$ occurs in $\pi_z$, so $\pi_z-\pi_x$ is
  a character. Thus (b) follows from Lemma \ref{L:EF}(a), and (c)
  follows from applying Lemma \ref{L:EF}(b) to an inertia generator of
  a place of $L$ above $t\mapsto\infty$.
\end{proof}

\begin{Remark} A weakening of doubly transitivity is primitivity. It
  is easy to see that an analog does not hold for primitive Galois
  groups. For instance let $5\le m$ and $2\le k\le m-2$ be integers
  with $2k\ne m$. Set $g(Z)=Z^m-Z$. Then
  $A:=\Gal(g(Z)-t/\QQQ(t))=\sym_m$.  Let $L$ be a splitting field of
  $g(Z)-t$ over $\QQQ(t)$, and $A_x\isom\sym_k\times\sym_{m-k}$ be a
  setwise stabilizer of a $k$-set in $\{1,2,\dots,m\}$. Let $f(t,X)$
  be a minimal polynomial of a primitive element of the fixed field of
  $A_x$ over $\QQQ(t)$.  From this setting we obtain that $f(g(Z),X)$
  is reducible in $\QQQ[Z,X]$. Therefore $\Red_f(\ZZZ)$ is an infinite
  set. Furthermore, the genus of the curve $f(T,X)=0$ goes to infinity
  with increasing $m$. For instance if $m$ is prime, then this genus
  is $1+\binom{m}{k}\frac{mk-k^2-m-1}{2m}\ge1$.
\end{Remark}

\begin{Remark} It does not seem to be obvious that we can replace the
  conclusion (b) in Theorem \ref{T:doubly}, namely that $k(t,x)$ has
  genus $0$, by the stronger conclusion that $k(t,x)$ is rational. An
  attempt to prove this stronger property leads to an interesting
  arithmetic question: Suppose that the assumptions of Theorem
  \ref{T:doubly} hold, but that $k(t,x)$ is not rational. We use the
  notation from the proof of Theorem \ref{T:doubly}. There we have
  seen that $\pi_z-\pi_x$ is a character, so in particular
  $[k(z):k(t)]\ge [k(t,x):k(t)]$. Let $\frakp_\infty$ be the rational
  place $t\mapsto\infty$. If $\frakp_\infty$ is totally ramified in
  $k(z)$, then so is this place in $k(t,x)$ by Lemma \ref{L:EF}(b), so
  the field $k(t,x)$ has a rational place, hence is rational. Thus
  there are two places of $\bk(z)$ above $\frakp_\infty$. Let $r$ and
  $s$ their ramification indices. By the argument above, there are two
  places of $\bk(t,x)$ above $\frakp_\infty$. As they are not
  rational, they are algebraically conjugate, so they have the same
  ramification index $u$. Clearly $r$ and $s$ divide $u$. The field
  degree estimation from above however gives $r+s\ge2u$. Thus $r=s=u$.
  We obtain $\pi_x=\pi_z$. Fields with this equality of permutation
  characters are said to be \emph{arithmetically equivalent}, see
  \cite{Klingen:Buch} for a book devoted to this subject. Thus we are
  lead to the following

\begin{Question} Let $k$ be a field, and $L/k(t)$ a finite Galois
  extension of the rational field $k(t)$. Let $k(t)\le k(z)\le L$ be a
  rational field, and $k(t)\le E\le L$ be a field which is
  arithmetically equivalent to $k(z)$ over $k(t)$. Does this imply
  that $E$ is a rational field as well?
\end{Question}

\end{Remark}

\subsection{Rational Specializations} An essential tool in our
investigation of integral Hilbert sets is Siegel's theorem about
algebraic curves with infinitely many integral points, combined with
the ramification behavior above infinity of Siegel functions. If we
look at rational specializations, then an analog of Proposition
\ref{P:Hset} holds, where Siegel's theorem is replaced by Falting's
theorem that a curve with infinitely many $k$-rational points has
genus at most $1$. The proof of the following proposition is similar
to the proof of Proposition \ref{P:Hset}, but simpler because we need
not worry about integrality.

\begin{Proposition} Let $k$ be a field which is finitely generated
  over $\QQQ$. Let $f(t,X)\in k(t)[X]$ be irreducible. Suppose that
  $f(\bar t,X)$ is reducible for infinitely many $\bar t\in k$. Then
  the splitting field $L$ of $f(t,X)$ over $k(t)$ contains a field
  $E\supset k(t)$ such that
  \begin{itemize}
  \item[(a)] $f(t,X)$ is reducible over $E$.
  \item[(b)] $E$ is either a rational field, or the function field of
    an elliptic curve with positive Mordell-Weil rank.
  \end{itemize}
\end{Proposition}

An application, whose proof is completely analogously to the proof of
Theorem \ref{T:doubly}, is the following finiteness statement.

\begin{Theorem} Let $k$ be a field which is finitely generated
  over $\QQQ$. Let $f(t,X)\in k(t)[X]$ be irreducible, with Galois
  group acting doubly transitively on the roots of $f(t,X)$.  If
  $f(\bar t,X)$ is reducible for infinitely many $\bar t\in k$, then
  the following holds, where $x$ is a root of $f(t,X)$:
\begin{enumerate}
\item[(a)] $f(t,X)$ is absolutely irreducible, and
\item[(b)] $k(t,x)$ has genus $\le1$.
\end{enumerate}
\end{Theorem}

\subsection{The prime degree case over the rationals}\label{SS:PrimeDeg}

Here we look at the case that the degree of the irreducible polynomial
$f(t,X)$ in $X$ is a prime number $p$. Let $A$ be the Galois group of
$f(t,X)$. By a classical result of Burnside (see e.g.\ \cite[Theorem
XII.10.8]{Huppert3}, \cite[Theorem 3.5B]{DixMort}), either $C_p\le
A\le\AGL_1(p)$, a case immediately dealt with, or $A$ is doubly
transitive. Though we treated the doubly transitive Galois groups in
the previous section, there are a few more things we can do in the
prime degree case.

The theorem below is an extension of \cite[Theorem 1.2]{PM:HIT}, the
method is different though.

Note that $h(X')-t$ with $h(X')\in \ZZZ[X']$ has a root in $\QQQ$ for
each $\bar t\in h(\ZZZ)$. If $x'$ is a root of $h(X')-t$, and $x\in
\QQQ(x')$ with $\QQQ(t,x)=\QQQ(x')$, then the minimal polynomial
$f(t,X)$ of $x$ over $\QQQ(t)$ has a root for the same (up to finitely
many exceptions) specializations $\bar t\in h(\ZZZ)$. The following
result shows that the converse holds in the odd prime degree case.

\begin{Theorem}\label{T:p} Let $f(t,X)\in \QQQ(t)[X]$ be irreducible of
  prime degree $p\ge3$ in $X$. Suppose that $f(\bar t,X)$ is reducible
  for infinitely many $\bar t\in\ZZZ$. Let $x$ be a root of $f(t,X)$.
  Then there is $x'\in\QQQ(t,x)$, such that $\QQQ(t,x)=\QQQ(x')$ and
  $t=h(x')$ with $h(X')\in\QQQ[X']$.
\end{Theorem}

\begin{proof} We use Lemma \ref{L:RedGal}. An element in
  $A$ of order $p$ is a transitive $p$-cycle on $A/A_x$. But $A_z$ is
  intransitive on $A/A_x$, so the order of $A_z$ is not divisible by
  $p$. But $p=[A:A_x]$, so $p$ must divide $[A:A_z]$. Let $\sigma$ be
  a generator of the inertia group $I$. So $\sigma$ has $m$ cycles of
  equal length on $A/A_z$, with $m=1$ or $2$. As $p$ is odd, $p$
  divides these cycle lengths, so in particular $\sigma$ has order
  divisible by $p$ on $A/A_z$. Thus $\sigma$ acts as a $p$-cycle on
  $A/A_x$. This means that the rational place $t\mapsto\infty$ is
  totally ramified in $\QQQ(t,x)$. If $A$ is doubly transitive, then
  $\QQQ(t,x)$ in addition has genus $0$ by Lemma \ref{L:EF}, and is a
  rational field because the unique place above $t\mapsto\infty$ must
  be rational, the claim follows in this case.
  
  Thus suppose that $A$ is not doubly transitive. Then $C_p\le
  A\le\AGL_1(p)$ in its action on $A/A_z$. An intransitive subgroup of
  such a group fixes a point, therefore $A_x$ is contained in a
  conjugate of $A_z$. So the fixed field $\QQQ(t,x)$ of $A_x$ has
  again genus $0$, and we complete the argument as above.
\end{proof}

\begin{Remark*} The above proof fails for $p=2$, because we cannot
  conclude that $p$ divides the cycles lengths of $\sigma$ on
  $A/A_z$. Indeed, the theorem does not hold for $p=2$. A
  counterexample is $f(t,X)=X^2+X-dt^2$ for a squarefree integer
  $d>1$.
\end{Remark*}

\subsection{Primitive Galois Groups}\label{SS:Primitive}

While we got a reasonably smooth result about Hilbert sets of
polynomials with a doubly transitive Galois group, the results are
less pleasant if we weaken the assumption on the Galois group to be
merely primitive. This section contains those results which we
achieved without using the classification of the finite simple groups.

\begin{Definition} Let $k$ be a finitely generated field of
  characteristic $0$. Denote by $\text{CF}(k)$ the set of those
  non-abelian simple groups which appear as composition factors of
  $\Gal(g(Z)-t/k(t))$ for Siegel functions $g(Z)$ over $k$.
\end{Definition}

In a bigger project \cite{PM:MonSiegel}, the simple groups
classification has been used to determine the sets $\text{CF}(k)$. In
particular, we obtained that, except for the alternating groups,
$\text{CF}(k)$ is finite. We come back to this
in Section \ref{S:AppClassif}.

\begin{Theorem}\label{T:main}
  Let $k$ be a finitely generated field extension of $\QQQ$, and $R$ a
  finitely generated subring of $k$. Let $f(t,X)\in k(t)[X]$ be
  irreducible, and assume that the Galois group $A$ of $f(t,X)$ over
  $k(t)$ acts primitively on the roots of $f(t,X)$. Suppose
  furthermore that $A$ has a non-abelian composition factor which is
  not contained in $\text{CF}(k)$. Then $\Red_f(R)$ is finite.
\end{Theorem}

The following result is, in terms of composition factors, a converse
to the previous theorem.

\begin{Theorem}\label{T:mainU}
  Let $k$ be a finitely generated field extension of $\QQQ$. Let
  $\Si\in \text{CF}(k)$ and $a\in\NNN$ be arbitrary. Then there exist an
  irreducible polynomial $f(t,X)\in k(t)[X]$ and a finitely generated
  subring $R$ of $k$, such that the following holds:
\begin{itemize}
\item[(a)] $\abs{\Red_f(R)}=\infty$.
\item[(b)] $\Gal(f(t,X)/k(t))$ acts primitively on the roots of
  $f(t,X)$.
\item[(c)] $\Si$ is a composition factor of $\Gal(f(t,X)/k(t))$.
\item[(d)] The genus of the curve $f(T,X)=0$ is $>a$.
\end{itemize}
\end{Theorem}

\begin{proof}[Proof of Theorem \ref{T:main}] We use Lemma
  \ref{L:RedGal}. We merely need to show that $A$ acts faithfully on
  $A/A_z$. Suppose the action is not faithful. Then there is a
  non-trivial normal subgroup $N\triangleleft A$ with $N\le A_z$. By
  primitive and faithful action of $A$ on $A/A_x$ we get $A=NA_x$.
  However, Lemma \ref{L:RedGal}(a) says that $A_zA_x$ is a proper
  subset of $A$. But $A_zA_x=A_z(NA_x)=A_zA=A$, a contradiction.
\end{proof}

In order to prove Theorem \ref{T:mainU} we need an easy group
theoretic observation.

\begin{Lemma}\label{L:wreathnormalizer} Let $G$ be a primitive
  non-regular permutation group on a set $\Delta$. Let $p$ be a prime,
  and $C_p\le H\le\AGL_1(p)$. Let $W=G^p\rtimes H$ be the wreath
  product, in the natural imprimitive action on the disjoint union of
  $p$ copies of $\Delta$. Let $\tilde W$ be a group acting on the same
  points, and suppose that $W$ is a normal subgroup of $\tilde W$.
  Then $\tilde W$ acts imprimitively, respecting the given system of
  imprimitivity of $W$. Therefore $\tilde W$ acts naturally and
  transitively on the cartesian product $\Delta^p$. This action of
  $\tilde W$ is primitive.\end{Lemma}

\begin{proof} $K=G^p$ is the kernel of the action of $W$ on the system
  of imprimitivity. Suppose there is $a\in \tilde W$ with $K\ne K^a$. We
  distinguish two cases. First suppose $K\cap K^a\ne1$. Then, by
  primitivity of $G$, $K\cap K^a$ is transitive on at least one and
  hence on each block $\Delta$. In particular, the orbits of $K^a$ are
  unions of $K$-orbits. On the other hand, $K$ and $K^a$ have the same
  number of orbits, so the blocks $\Delta$ are exactly the $K^a$
  orbits. Thus $K=K^a$, a contradiction.
  
  Next assume that $K\cap K^a=1$. Then $K^a$ acts faithfully and
  transitively as a normal subgroup of $\AGL_1(p)$ on the system of
  imprimitivity. So $p$ divides the order of $K^a\isom G^p$, but $p^2$
  does not. This contradiction shows that $K$ is normal in $\tilde W$.
  
  As the blocks $\Delta$ are the $K$-orbits, we obtain that $\tilde W$
  respects that system. By \cite[Lemma 2.7A]{DixMort}), the action of
  $W$ on $\Delta^p$ is primitive, so this is even more true for
  $\tilde W$.\end{proof}

In general, the composition of Siegel functions is not a Siegel
function. Conversely, if we write a Siegel function as a composition
of rational functions, then not all these rational functions need to
be Siegel functions. The following lemma clarifies this issue.

\begin{Lemma}\label{L:SiegelComp} Let $k$ be a finitely generated
  field over $\QQQ$, and $g(Z)\in k(Z)$ a Siegel function over $k$ of
  degree $>1$. Then there is a decomposition $g(Z)=a(b(Z))$ with
  $a,b\in k(Z)$, such that the following holds:\begin{itemize}
  \item[(a)] $a(Z)$ is a functionally indecomposable Siegel function
    of degree $>1$.
  \item[(b)] There is $0\ne\delta\in k$, such that $\delta b(Z)$ is a
    Siegel function, or the followings holds: There are linear
    fractional functions $\lambda,\mu\in k_1(Z)$ over a quadratic
    extension $k_1$ of $k$, $m\in\NNN$, with $b(Z)=\lambda(\mu(Z)^m)$.
    In this case, $\Gal(b(Z)-t/k(t))$ is solvable.
\end{itemize}
\end{Lemma}

\begin{proof} Let $R$ be a finitely generated subring of $k$ such that 
  $\abs{g(k)\cap R}=\infty$.
  
  Write $g(Z)=a(b(Z))$ with $a,b\in k(Z)$ and $a(Z)$ being
  indecomposable over $k$ of degree $>1$. As $b(k)\subseteq
  k\cup\{\infty\}$, we have $\abs{a(k)\cap R}=\infty$, so (a) clearly
  holds.
  
  From $\abs{g^{-1}(\infty)}\le2$ we obtain
  $\abs{a^{-1}(\infty)}\le2$.
  
  We first analyze the case $\abs{g^{-1}(\infty)}=1$. Because
  $\Gal(\bar k/k)$ acts on $g^{-1}(\infty)$, this single element in
  this fiber must be rational or $\infty$. By a linear fractional
  change we may assume that $g^{-1}(\infty)=\{\infty\}$, so $a(Z)$ is
  a polynomial. First suppose $k\ne\QQQ$. Write
  $a(Z)=a_rZ^r+s_{r-1}Z^{r-1}+\dots+a_1Z+a_0$ with $a_i\in k$. By
  assumption, there are infinitely many $\bar z\in k$ such that
  $\beta=b(\bar z)$ fulfills $a(\beta)=\rho\in R$. So $\beta$ is
  integral over
  $R'=R[\frac{a_n-1}{a_n},\dots,\frac{a_1}{a_n},\frac{a_0}{a_n}]$.
  Thus replace $R$ by a finitely generated ring containing the
  integral closure of $R'$ in $k$, using \cite[Chapter 2, Proposition
  4.1]{Lang:Dio}. So $b(Z)$ is a Siegel function with respect to this
  ring. Thus assume $k=\QQQ$, so $a_i\in\QQQ$. Let $w$ be a common
  multiple of the denominators of $a_i$. The previous consideration
  shows that $a_nw\beta\in\QQQ$ is integral over $\ZZZ$, hence
  contained in $\ZZZ$. So $a_nwb(Z)$ assumes infinitely many integral
  values on $\QQQ$, and (a) follows again.
  
  Now assume $\abs{g^{-1}(\infty)}=2$. Write
  $g^{-1}(\infty)=\{\lambda_1,\lambda_2\}$. Then the $\lambda_i$ are
  either in $k\cup\{\infty\}$, or they generate a quadratic extension
  $k_1$ of $k$. Furthermore, we obtain $b^{-1}(\lambda_i)=\{\mu_i\}$,
  with $\mu_i\in k\cup\{\infty\}$ in the former case, or $\mu_i\in
  k_1$ in the latter case. At any rate, there are linear fractional
  functions $\lambda,\mu\in k_1(Z)$ such that, with $\tilde
  b(Z):=\lambda^{-1}(b(\mu^{-1}(Z)))$, the following holds: $\tilde
  b^{-1}(\infty)=\{\infty\}$, $\tilde b^{-1}(0)=\{0\}$, and $\tilde
  b(1)=1$. This implies $\tilde b(Z)=Z^m$. The Galois group of
  $b(Z)-t$ over $k_1(t)$ is the same one as the Galois group of
  $\tilde b(Z)-t$ over $k_1(t)$. This group is contained in
  $\AGL_1(n)$, hence solvable. The Galois group of $g(Z)-t$ over
  $k(t)$ is an extension of the former Galois group by at most the
  index $2$, so is solvable as well. This proves (b).
\end{proof}

\begin{Corollary}\label{C:Si} Let $k$ be a finitely generated
  field extension of $\QQQ$, and $g(Z)\in k(Z)$ a Siegel function over
  $k$. Let $\Si$ be a non-abelian composition factor of
  $\Gal(g(Z)-t/k(t))$. Then there is a functionally indecomposable
  Siegel function $\tilde g(Z)$ over $k$, such that $\Si$ is a
  composition factor $\Gal(\tilde g(Z)-t/k(t))$.
\end{Corollary}

\begin{proof} If $g(Z)=g_1(g_2(\dots g_r(Z)\dots))$ with functionally
  indecomposable rational functions $g_i(Z)\in k(Z)$, then $\Si$ is a
  composition factor of $\Gal(g_i(Z)-t/k(t))$ for some index $i$. See
  Glauberman's argument in \cite[Prop.~2.1]{GT} for this fact which is
  less obvious than it might appear at a first glance.

The assertion now follows from Lemma \ref{L:SiegelComp}
\end{proof}

\begin{proof}[Proof of Theorem \ref{T:mainU}] Let $\Si\in
  \text{CF}(k)$, so $\Si$ is a non-abelian composition factor of
  $\Gal(g(Z)-t/k(t))$ for a Siegel function $g(Z)$ over $k$. By
  Corollary \ref{C:Si} we may assume that $g$ is functionally
  indecomposable. The Galois group $G$ of $g(Z)-t$ over $\bar k(t)$
  acts primitively on the roots of $g(Z)-t$, because $g(Z)$ is
  functionally indecomposable over $\bar k$ by Theorem \ref{T:absind}.
  Let $R$ be a finitely generated ring in $k$ with $\abs{g(k)\cap
    R}=\infty$.
  
  Let $p$ be a prime. Choose $\alpha\in R$ such that the following
  holds: $0$ and $\infty$ are not branch points of $g(Z)-\alpha$, and
  the $p$-th powers of the branch points of $g(Z)-\alpha$ are all
  distinct. These general position assumptions will be used in a genus
  computation below. Set $\tilde g(Z):=(g(Z)-\alpha)^p$. Let $\zeta$
  be a primitive $p$-th root of unity. By our choices, the sets of
  branch points of the splitting fields of the $p$ functions
  $g(Z)-\alpha-\zeta^it^{1/p}$, $i=1,2,\dots,p$ over $\bar k(t^{1/p})$
  are pairwise disjoint. As $\bar k(t^{1/p})$ does not posses
  unramified finite extensions, each of these splitting fields is
  linearly disjoint to the compositum of the remaining $p-1$ ones. This
  implies that the Galois group $W$ of $\tilde g(Z)-t$ over $\bar
  k(t)$ is the wreath product $G^p\rtimes C_p$.
  
  Let $L$ be a splitting field of $\tilde g(Z)-t$ over $k(t)$, and
  $\hat k$ the algebraic closure of $k$ in $L$. Again $W=\Gal(L/\hat
  k(t))$. Set $\tilde W=\Gal(L/k(t))$. By Lemma
  \ref{L:wreathnormalizer}, $\tilde W$ has a maximal subgroup $V$,
  such that $V$ is intransitive on the roots of $\tilde g(Z)-t$.
  Indeed, one orbit of $V$ has length $p$.
  
  We have $\abs{\tilde g(k)\cap R}=\infty$. Let $f(t,X)\in k(t)[X]$ be
  a minimal polynomial of a primitive element of the fixed field of
  $V$ in $L$ over $k(t)$. We have verified (a), (b), and (c) of our
  theorem.
  
  It remains to compute the genus of $f(t,X)=0$. Recall that $f(t,X)$
  is absolutely irreducible by Corollary \ref{C:HIT:absirr}. We work
  over $\hat k(t)$. Let $G_1$ be the stabilizer of a point in the
  given action of $G$ on the roots of $g(Z)-t$. Let $x$ be a root of
  $f(t,X)$. The stabilizer $V\cap W$ in $W$ of $x$ can be identified
  with $W_1:=G_1^p\rtimes C_p$. Note that if $n$ is the degree of
  $g(Z)$, then $f(t,X)$ has degree $n^p$.
  
  We take advantage of the general position assumptions of the
  branching locus of $\hat g(Z)-t$ in order to get an easy genus
  computation using the Riemann--Hurwitz genus formula. Let $\sigma_1$
  and $\sigma_2$ be inertia generators belonging to $t\mapsto0$ and
  $t\mapsto\infty$, and let $\tau_1,\dots,\tau_r$ be inertia
  generators coming from the branch points of $g(Z)$. Let $\ind$ refer
  to the action on $W/W_1$. Then $\sigma_i$ has precisely $n$ fixed
  points, and moves the remaining $n^p-n$ points in $p$--cycles. Thus
  $\ind(\sigma_i)=(n^p-n)(1-1/p)$. If the inertia generator belonging
  to $\tau_i$ has orbit lengths $v_1,v_2,\dots,v_s$ on the roots of
  $g(Z)-t$, then $\tau_i$ has the same orbit lengths on $W/W_1$, but
  each one occurs $n^{p-1}$ times. As $\hat k(Z)/\hat k(g(Z))$ is an
  extension of genus $0$ fields, we obtain
\[
\sum_{i=1}^r\ind(\tau_i)=n^{p-1}(2(n-1)).
\]
If $g_f$ is the genus of $\hat k(t,x)$, then
\begin{align*}
2(n^p-1+g_x) &= \ind(\sigma_1)+\ind(\sigma_2)+\sum_{i=1}^r\ind(\tau_i)\\
             &= 2(n^p-n)(1-\frac{1}{p})+n^{p-1}(2(n-1)),
\end{align*}
so
\[
g_x = (n^{p-1}-1)\frac{np-n-p}{p} > 0,
\]
and clearly $g_x\to\infty$ for $p\to\infty$.
\end{proof}

\section{Applying the simple groups classification}\label{S:AppClassif}

So far we have used only easy arithmetic and geometric properties of
the ramification structure of Siegel functions. In particular, the
results so far are not based on the classification of the finite
simple groups. In order to obtain more results, it is indispensable to
obtain good information on the Galois groups of $g(Z)-t$ for Siegel
functions $g$. This has been carried out in \cite{PM:MonSiegel}. There
we classify the possible Galois groups, and study which cases live
over the rationals. We quote three corollaries from this
classification.

\begin{Theorem} Let $k$ be a field of
  characteristic $0$, and $g(Z)\in k(Z)$ a Siegel function. Then each
  non-abelian composition factor of $\Gal(g(Z)-t/k(t))$ is isomorphic
  to one of the following groups: $\alt_j$ ($j\ge5$), $\PSL_2(7)$,
  $\PSL_2(8)$, $\PSL_2(11)$, $\PSL_2(13)$, $\PSL_3(3)$, $\PSL_3(4)$,
  $\PSL_4(3)$, $\PSL_5(2)$, $\PSL_6(2)$, $\M11$, $\M12$, $\M22$,
  $\M23$, $\M24$.
\end{Theorem}

\begin{Theorem}\label{T:MonSiegelQ} Let
  $g(Z)\in\QQQ(Z)$ be a Siegel function over $\QQQ$. Then each
  non-abelian composition factor of $\Gal(g(Z)-t/\QQQ(t))$ is
  isomorphic to one of the following groups: $\alt_j$ ($j\ge5$),
  $\PSL_2(7)$, $\PSL_2(8)$.
\end{Theorem}

\begin{Theorem} Let $g(Z)\in\QQQ(Z)$ be a
  Siegel function over $\QQQ$. Assume that $A=\Gal(g(Z)-t/\QQQ(t))$ is
  a simple group. Then $A$ is isomorphic to an alternating group or
  $C_2$.
\end{Theorem}

An immediate application of the latter theorem and Lemma
\ref{L:RedGal} is

\begin{Corollary}\label{C:mainQ:Galois_2}
  Let $f(t,X)\in\QQQ(t)[X]$ be irreducible with Galois group $A$,
  where $A$ is a simple group not isomorphic to an alternating group
  or $C_2$. Then $\Gal(f(\bar t,X)/\QQQ)=G$ for all but finitely many
  specializations $\bar t\in\ZZZ$.
\end{Corollary}

\begin{Remark*} This corollary becomes completely wrong if we allow
  rational specializations $\bar t\in\QQQ$. Indeed, many interesting
  simple groups are Galois groups of polynomials $A(X)-tB(X)$ with
  $A,B\in\QQQ[X]$, see \cite[Appendix, Table 10]{MM}. So for each
  specialization $\bar t=A(\bar z)/B(\bar z)$ with $\bar z\in\QQQ$ we
  obtain a smaller Galois group, because $A(X)-\bar tB(X)$ becomes
  reducible.
\end{Remark*}

Similarly, Theorems \ref{T:main} and \ref{T:MonSiegelQ} give

\begin{Corollary}
  Let $f(t,X)\in\QQQ(t)[X]$ be irreducible, and assume that the Galois
  group of $f(t,X)$ over $\QQQ(t)$ acts primitively on the roots of
  $f(t,X)$ and has a non-abelian composition factor which is not
  alternating and not isomorphic to $\PSL_2(7)$ or $\PSL_2(8)$. Then
  $\Red_f(\ZZZ)$ is finite.
\end{Corollary}

\noindent{\sc IWR, Universit\"at Heidelberg, Im Neuenheimer Feld 368,\\
69120 Heidelberg, Germany}\par
\noindent{\sl E-mail: }{\tt Peter.Mueller@iwr.uni-heidelberg.de}
\end{document}